\documentclass[english]{amsart}
\usepackage{a4wide}

\usepackage{ulem}
\usepackage{graphicx}

\usepackage{amsmath}

\usepackage{amssymb}

    \DeclareFontFamily{U}{wncy}{}
    \DeclareFontShape{U}{wncy}{m}{n}{<->wncyr10}{}
    \DeclareSymbolFont{mcy}{U}{wncy}{m}{n}
    \DeclareMathSymbol{\Sha}{\mathord}{mcy}{"58}

\usepackage{amscd}
\usepackage[all,cmtip,line]{xy}

\newlength{\ownl}

\newcommand{\Spec}{{\operatorname{Spec}\,}}

\newcommand{\R}{{\mathbb{R}}}

\newcommand{\CO}{{\mathcal{O}}}

\newcommand{\cF}{\mathcal{F}}

\newcommand{\cO}{\mathcal{O}}

\def\RCS$#1: #2 ${\expandafter\def\csname RCS#1\endcsname{#2}}
\RCS$Revision: 327 $
\RCS$Date: 2010-11-15 13:46:16 -0600 (Mon, 15 Nov 2010) $

\usepackage{amsthm}

\newtheorem{thm}{Theorem}[subsection]

 \newtheorem{lemma}[thm]{Lemma}

 \theoremstyle{definition}
 \theoremstyle{definition}
 \theoremstyle{remark}

\numberwithin{equation}{subsection}

\theoremstyle{definition}

\usepackage[bookmarksopen,bookmarksdepth=2,breaklinks=true]{hyperref}

\begin{document}
\title[]{Schemes in Lean}

\author{Kevin Buzzard} \email{k.buzzard@imperial.ac.uk} \address{Department of
  Mathematics, Imperial College London}
\author{Chris Hughes} \email{christopher.hughes17@imperial.ac.uk} \address{Department of
  Mathematics, Imperial College London}
\author{Kenny Lau} \email{kin.lau17@imperial.ac.uk} \address{Department of
  Mathematics, Imperial College London}
\author{Amelia Livingston} \email{amelia.livingston17@imperial.ac.uk} \address{Department of
  Mathematics, Imperial College London}
\author{Ramon Fern\'andez Mir} \email{Ramon.FernandezMir@ed.ac.uk}
\address{School of Informatics, University of Edinburgh}
\author{Scott Morrison} \email{scott@tqft.net}
\address{School of Mathematics and Statistics, University of Sydney}

\thanks{The first author was supported in part by EPSRC grant EP/L025485/1.}

\begin{abstract}
  We tell the story of how schemes were formalised in three different ways in the Lean theorem prover.
\end{abstract}

\maketitle
\section{Introduction and overview}

\subsection{Varieties and schemes in algebraic geometry}

Before 1960, algebraic geometry was done via the theory of algebraic varieties, finite-dimensional objects defined over a fixed algebraically closed field, or ``universal domain''. The standard reference text was Weil's 1946 book ``Foundations of algebraic geometry'' ~\cite{weil}, and in the final chapter ``Comments and discussions'', Weil remarks that ``it would be very convenient to have\ldots a principle of reduction modulo~$p$ '', a phenomenon which Weil would have known well should exist but which was extremely inconvenient to do in this setting.

Schemes were introduced by Grothendieck in the 1960s (following earlier ideas of Chevalley, and building on ideas of Zariski) as the building blocks for a new algebraic geometry. Grothendieck did not need to work over a fixed base field; his foundations worked with general commutative rings rather than algebraically closed fields, enabling reduction modulo~$p$ to become possible. Given a general ``geometric object'' (for example, a topological space), one can consider the set of continuous real-valued functions on this object, and pointwise addition and multiplication turn this space of functions into a commutative ring. Grothendieck's observation was that one could make a construction in the opposite direction: starting with an arbitrary commutative ring~$R$, he constructed a geometric object $\Spec(R)$ called an \emph{affine scheme}; this was a topological space equipped with a sheaf of functions $\cO_X$ on this space, such that $\cO_X(\Spec(R))$, the ``allowable'' functions on~$\Spec(R)$, was~$R$ again. Grothendieck defined a general \emph{scheme} by gluing affine schemes together (following Weil, who in~\cite{weil} had defined an abstract variety by gluing affine varieties together\footnote{Note that Whitney had defined an abstract real manifold ten years earlier in~\cite{whitney}, by gluing together open subsets of Euclidean space in the now familiar way.}). This new viewpoint proved incisive -- ten years later the Weil conjectures, questions about the number of points which algebraic varieties have over finite fields, were being proved using machinery such as \'etale cohomology, which had been developed using schemes.

\subsection{Formalising schemes -- the history}

As mentioned in the abstract, schemes were formalised three times in Lean \cite{DeMouraKongAvigadVanDoornvonRaumer}; each formalisation was better-behaved than the one before.

The first formalisation was by the first three authors (KB, CH, KL), and it evolved in late 2017 and early 2018 from the Thursday evening Xena Project meetings at Imperial College London, where undergraduates learn how to formalise mathematics in Lean (KB is the staff member running the meetings; CH and KL were, at the time, first year mathematics undergraduates). It started with KL formalising the theory of localisation of rings as a project, and with KB suggesting that schemes would be a natural way to take the theory further in Lean. It was their first attempt to formalise anything non-trivial in a theorem prover (and its ultimate success inspired KB to see exactly how far Lean could be pushed, eventually resulting in~\cite{perfectoid}). However some poor (in retrospect) design decisions were made when it came to the theory of localisation, and in this first iteration these decisions resulted in messy infrastructure which would not scale. It quickly became clear that an extensive rewrite was needed. It also became apparent that such a formalisation had never been embarked upon before in any other theorem prover (parts of the theory of affine schemes had been formalised, but nothing approaching the definition of a scheme) -- this perhaps says something about the interests of mathematical formalisation community, at least pre-2017. The (abandoned) project is still currently online at~\cite{schemes-v1}.

AL and RFM then joined the project, with AL developing a very robust theory of localisation in the form which was needed, and RFM rewriting the definition of a scheme from scratch as part of his 2018--2019 MSc project supervised by KB. RFM ultimately produced a definition which was usable -- indeed several other basic results from the Stacks project website~\cite{stacks-project} and Hartshorne's algebraic geometry textbook~\cite{hartshorne} were also proved in this iteration, mostly by KL. This version is currently online at~\cite{schemes-v2}.

However, sheaves were defined ``by hand'' in this iteration; there was one definition for a sheaf of types (or sheaf of sets, depending on your foundations), one for a sheaf of abelian groups and another one for a sheaf of rings. The contribution of SM was to build enough abstract category theory in Lean to enable a third definition using this category-theoretic language. At this point a design change was also made; the definition of the sheaf on $\Spec(R)$ was changed from the definition in~\cite{stacks-project} to the equivalent definition in~\cite{hartshorne}. This was the version which finally made it into Lean's mathematics library, in commit {\tt b79fc0379ae786153fc22ce5ee6751505e36a3d9} of {\tt mathlib}~\cite[{\href{https://github.com/leanprover-community/mathlib/blob/master/src/algebraic_geometry/Scheme.lean}{\tt src/algebraic\_geometry/Scheme.lean}}]{mathlib}. The current {\tt mathlib} documentation for schemes (perhaps more readable for people who do not want to look directly at the code) is available \href{https://leanprover-community.github.io/mathlib_docs/algebraic_geometry/Scheme.html}{at the Lean community website}.

\subsection{Organisation of the paper.}
Most of the formalisation was plain sailing, but occasionally we ran into unexpected problems. The layout of this paper is as follows. In the next section we go over the details of some of the mathematics involved. In the three sections after that we explain the three approaches to formalising the material, emphasising the parts which did not go smoothly and explaining how this affected future design decisions.

\section{Mathematical details.}

Convention: all rings are commutative and have a 1; all ring homomorphisms send~1 to~1.

Let~$X$ be a topological space. A \emph{presheaf of rings} $\cF$ on~$X$ is a way to associate a ring $\cF(U)$ to each open subset $U\subseteq X$, and a ring homomorphism $\rho_{UV}:\cF(U)\to\cF(V)$ to each inclusion $V\subseteq U$ of open subsets of~$X$, such that $\rho_{UU}$ is the identity for all opens $U$, and $\rho_{WU}=\rho_{VU}\circ\rho_{WV}$ if $U\subseteq V\subseteq W$. In other words, $\cF$ is a contravariant functor from the category of open subsets of~$X$ to the category of rings. If $f\in\cF(U)$ and $V\subseteq U$, we write $f|_V$ as shorthand for $\rho_{UV}(f)$. The model example to keep in mind is where $\cF(U)$ is defined to be the ring of continuous functions $U\to\R$, and $\rho_{VU}$ sends a continuous function on~$U$ to its restriction to~$V$.

A presheaf of rings~$\cF$ is said to be a \emph{sheaf of rings} if elements of $\cF(U)$ can ``be defined locally''. More precisely, a presheaf of rings $\cF$ is a sheaf of rings if for every open set $U$ and every open cover $U=\bigcup_{i\in I}U_i$ of $U$ by open sets $U_i$, the sequence
$$0\to\cF(U)\to\prod_{i\in I}\cF(U_i)\to\prod_{i,j\in I}\cF(U_i\cap U_j)\eqno{(\ddag)}$$
is exact, where the map $\cF(U)\to\prod_i\cF(U_i)$ sends $f$ to $(f|_{U_i})$, and the map $\prod_i\cF(U_i)\to\prod_{i,j}\cF(U_i\cap U_j)$ sends $(f_i)$ to the element whose value at $(i,j)$ is $f_i|_{U_i\cap U_j}-f_j|_{U_i\cap U_j}$. In words, the exactness of $(\ddag)$ is the assertion that if you have a collection of elements $f_i\in\cF(U_i)$ which agree on all overlaps $U_i\cap U_j$, then there is a unique $f\in\cF(U)$ whose restriction to $U_i$ is $f_i$ for all $i$. A topological space~$X$ equipped with a sheaf of rings~$\cO_X$ is called a \emph{ringed space}.

If $(X,\cO_X)$ is a ringed space, then to each point $u\in X$ one can associate the \emph{stalk} $\cO_{X,x}:=\varinjlim_{x\in U}\cO_X(U)$ of ``functions defined near~$x$'', a filtered colimit of rings and hence also a ring. If this ring is local (that is, has a unique maximal ideal) for every $x\in X$ we say that $(X,\cO_X)$ is a \emph{locally ringed space}.

A fundamental construction in this area is the following. Given a ring~$R$, let $\Spec(R)$ denote the set of its prime ideals. The Zariski topology is a natural topology on $\Spec(R)$: for $I$ an ideal we define $V(I)$ to be the prime ideals of~$R$ containing~$I$, and the $V(I)$ as~$I$ varies are the closed sets of $\Spec(R)$. We shall explain later how to use the theory of localisations of rings to construct a presheaf of rings $\CO_R$ on $\Spec(R)$ and to prove that it is a sheaf. The ringed space $(\Spec(R), \cO_R)$ is locally ringed, and is called the \emph{affine scheme} associated to the ring~$R$.

A \emph{scheme} is a locally ringed space $(X,\cO_X)$ with the property that $X$ can be covered by opens $U_i$ such that the induced locally ringed space $(U_i,\cO_X|_{U_i})$ is isomorphic to $(\Spec(R_i),\cO_{R_i})$ for some ring $R_i$.

If $R$ is a ring then $\Spec(R)$ is a scheme, because it can be covered by the open affine subset $U=\Spec(R)$.

\medskip

To KB's naive eyes in 2017, everything here looked straightforward. But in fact there were several gotchas involved in turning this into a formal definition. The rest of this paper describes them, and how they were dealt with in Lean.

\section{The first definition.}

We begin with a discussion of localisation, and the construction of the sheaf of rings on the topological space associated to a ring. Let $R$ be a ring and let $S$ be a submonoid of $(R,\times)$ (for some reason mathematicians often call $S$ a ``multiplicative subset''; it contains $1$ and is closed under multiplication so it is precisely a submonoid). The \emph{localisation} $R[1/S]$ of $R$ at $S$ is a ring equipped with a canonical map $i:R\to R[1/S]$ and having the following universal property: for any ring $A$ and any map $f:R\to A$ such that $f(s)\in A^\times$ for all $s\in S$, there is a unique extension $\tilde{f} : R[1/S]\to A$ of $f$ to $R[1/S]$ making the obvious diagram
$$\xymatrix
{
  R[1/S] \ar[dr]^{\exists!\tilde{f}}        &\\
  R\ar[u]^{i}\ar[r]^{f}&A
}$$
commute. The special case where $S$ is the submonoid $\{1,f,f^2,f^3,\ldots\}$ generated by $f\in R$ has its own notation $R[1/f]$.

Standard arguments involving universal objects show that $R[1/S]$ is uniquely defined up to unique isomorphism, if it exists. In fact, these latter arguments do not even assume that~$S$ is a submonoid; in general the localisation at a general subset $S$ of $R$ is the localisation at the submonoid of~$R$ generated by~$S$. Existence is shown by an explicit construction: one puts an appropriate equivalence relation on $R\times S$ (an element $(r,s)$ is thought of as representing the fraction $r/s$) and puts a ring structure on the quotient. All of this was formalised in Lean with little fuss.

The construction of the topological space $\Spec(R)$ was equally uneventful. Next, a decision had to be made about how to put a sheaf of rings on $\Spec(R)$; there are several constructions in the literature. The prevailing philosophy of KB at the time was ``always follow the Stacks project''.\footnote{The Stacks project~\cite{stacks-project} covers a huge swathe of algebraic geometry and its prerequisites. As the project was largely a collaboration by KB with undergraduates, the fact that the Stacks project is freely available online was helpful as it was accessible to everyone. Moreover, the Stacks project is careful to use `Tags', which never change, to refer to particular theorems or sections of the material, so it is easy to cross-reference.} The approach in the Stacks project is as follows. Let~$X$ be a topological space. In \cite[\href{https://stacks.math.columbia.edu/tag/009H}{Tag 009H}]{stacks-project} it is explained how to set up the theory of presheaves and sheaves on a basis for the topology on~$X$, and in \cite[\href{https://stacks.math.columbia.edu/tag/009N}{Tag 009N}]{stacks-project} it is shown how to extend a sheaf on a basis for~$X$ to a sheaf on~$X$. Formalisation of these results, and their extension to sheaves of rings, was straightforward. Indeed, this part of the formalisation went like a dream -- one line of very abstract categorical mathematics would correspond to one line of Lean code. These arguments reduced the question to defining a sheaf on the basis of $\Spec(R)$ consisting of open sets of the form $D(f)$, where $f\in R$ and $D(f)$ is the set of prime ideals not containing~$f$. A set of the form $D(f)$ is called a \emph{basic open set}. Note that $D(f)$ is open, because it is the complement of the closed set $V(I)$ where $I=(f)$.

\subsection{Definition of $\cO_X(D(f))$, and ``canonical maps''.}

This was a part of the story where things did not go as smoothly as envisaged. The first problem was this. Say $R$ is a ring, and $f,g\in R$ have the property that $D(f)=D(g)$, that is, a prime ideal contains~$f$ if and only if it contains~$g$. This can certainly happen in non-obvious ways -- for example if $x\in R$ and $f=x^2$ and $g=-x^3$, then a prime ideal contains $g$ if and only if it contains~$x$, and so if and only if it contains~$f$. The literature happily defines $\cO_X(D(f))$ to be $R[1/f]$ -- see for example \href{https://github.com/stacks/stacks-project/commit/453818f0a1fa56e0b830cec409289383ba1a38a1}{this commit} of \cite[\href{https://stacks.math.columbia.edu/tag/01HT}{Tag 01HT}]{stacks-project}:
\begin{quote}
If $f,g\in R$ are such that $D(f)=D(g)$, then by Lemma 01HS above there are canonical maps $M_f\to M_g$ and $M_g\to M_f$ which are mutually inverse. Hence we may choose any $f$ such that $U=D(f)$ and define $\tilde{M}(U)=M_f$.
\end{quote}
(Here $M_f$ is a module-theoretic generalisation of $R[1/f]$). No definition of ``canonical'' is supplied. The same is seen in Grothendieck's original work -- in~\cite{ega1}, section 1.3.1 of Chapter I, it is written
\begin{quote}
on sait que $A_f$ et $M_f$ s'identifient canoniquement \'a $S_f'^{-1}A$ et $S_f'^{-1}M$
\end{quote}
(again with no definition of ``canoniquement'') and in 1.3.3 we see ``un homomorphisme canonique fonctoriel $M_f\to M_g$'', and the first usage of ``$M_f=M_g$'' to denote a ``canonical'' isomorphism rather than a set-theoretic identity. Having established that $R[1/f]=R[1/g]$ (for this new ``canonically isomorphic'' meaning of $=$) Grothendieck defines $\cO_X(D(f))=R[1/f]$ and $\cO_X(D(g))=R[1/g]$ and continues.

This becomes an issue in Lean, because in our formalisation we have rings $R[1/f]$ and $R[1/g]$ which are most definitely \emph{not} equal; certainly the universal property gives isomorphisms between them, but one of these rings is a quotient of $R\times\{1,f,f^2,f^3,\ldots\}$ and the other is a quotient of $R\times\{1,g,g^2,g^3,\ldots\}$. If you like to think set-theoretically, you can think of an element of $R[1/f]$ as being a subset (an equivalence class) of $R\times\{1,f,f^2,\ldots\}$ and an element of $R[1/g]$ as being a subset of $R\times\{1,g,g^2,\ldots\}$. In particular, $R[1/f]$ and $R[1/g]$ are visibly not \emph{equal}.

One solution to this problem would be to throw the axiom of choice at it. For each open set $U\subseteq\Spec(R)$ which is known to be of the form $D(f)$ for some $f$ (possibly infinitely many), we choose some ``special'' $f_U$ such that $U=D(f_U)$ and define $\cO_X(U):=R[1/f_U]$. This of course works, but we envisaged that this definition would be frustrating to work with down the line.

Reid Barton on the Lean chat suggested the following approach instead, which was what we ultimately chose. If the open set $U$ is known to be of the form $D(f)$, then define $S_U$ to be the submonoid $\{g\in R\,|\,U\subseteq D(g)\}$ of elements of~$R$ which are non-zero on $U$, and define $\cO_X(U):=R[1/S_U]$. This has the advantage that no choices are involved; it bears the hallmark of a construction in constructive mathematics, where a rule of thumb is that if you can't figure out a natural way to choose one object from a set, then choose all of them.

This defines a presheaf of rings $\cO_X$ on the basic open subsets $D(f)$ of $\Spec(R)$, and by taking limits one can extend this definition to give a presheaf of rings on all of $\Spec(R)$. Amusingly, we then realised that to define schemes, this presheaf construction was all we needed.

\subsection{Defining schemes.}

We defined a scheme to be a topological space $X$ equipped with a sheaf of rings $\cO_X$, for which $X$ had a cover $X=\bigcup_iU_i$ such that each $U_i$ was isomorphic to $\Spec(R_i)$ for some ring $R_i$.

This isomorphism is in the books usually stated as an isomorphism of locally ringed spaces. However a locally ringed space is just a topological space equipped with a presheaf of rings and satisfying some extra axioms, so in particular 
an isomorphism of locally ringed spaces is the same as an isomorphism of spaces equipped with presheaves of rings. We have already defined a presheaf of rings on $\Spec(R)$, so we are almost done. It remains to define the pullback presheaf of rings on $U_i$, and this is easy: if $\iota: U_i\to X$ is the open immersion, then for $V\subseteq V_i$ open, the corresponding subset $\iota(V)\subseteq X$ is an open subset of~$X$, so we define the presheaf $\iota^*\cO_X$ by $\iota^*\cO_X(V)=\cO_X(\iota(V))$. It was straightforward to check that $\iota^*\cO_X$ is a presheaf of rings on~$U_i$, and our definition of a scheme was complete. Note that we did not need to demand that our original ringed space was locally ringed -- this follows from the fact that it is locally affine. In particular, our definition is mathematically equivalent to the usual definition, although the proof of this involves theorems which were not at this time formalised.

Our original definition was met with some scepticism by the computer scientists in the Lean community, however, and not for this reason above. A definition with no unit tests might contain an error. It was suggested that we prove a theorem about schemes, to provide evidence that our definition was correct. We decided to prove the theorem that an affine scheme was a scheme. Although the result sounds trivial, some work remains: a scheme is a space equipped with a sheaf of rings, and we have thus far only equipped $\Spec(R)$ with a presheaf of rings. To prove that affine schemes are schemes should be essentially equivalent to proving that the presheaf $\cO_X$ of rings on $\Spec(R)$ is a sheaf.

\subsection{Proving $\cO_X$ is a sheaf on $\Spec(R)$.}

Again, we followed the Stacks project. Here we ran into a serious technical issue, and one might even argue that this issue is typically overlooked in the literature -- it seems to be a very good example of a situation where mathematicians pay no attention to the details of what is happening, knowing that things are going to work out. When formalising, one has to check these details.

Let us explain the main problem. We need to check the sheaf axiom for $\cO_X$ considered as a sheaf on the basis $D(f)$ of open sets in $\Spec(R)$. A formal argument involving compactness reduces us to the case of a cover of a basic open by finitely many basic opens, so all that remains is checking that the sheaf axiom $(\ddag)$ holds for a cover of a basic open set by finitely many basic open sets.

The argument in \cite[\href{https://stacks.math.columbia.edu/tag/01HR}{Tag 01HR}]{stacks-project} now proceeds like this:

{\bf Step 1:} If $D(f)$ is a basic open in $\Spec(R)$, then $D(f)$ is canonically isomorphic to $\Spec(R[1/f])$, and the isomorphism identifies basic open subsets in $\Spec(R[1/f])$ with the basic open subsets in $\Spec(R)$ contained within $D(f)$. This reduces us to the case of checking that the sheaf axiom $(\ddag)$ holds for finite covers of $\Spec(R)$ by basic opens;

{\bf Step 2:} One can now translate the statement into a purely ring-theoretic lemma \cite[\href{https://stacks.math.columbia.edu/tag/00EJ}{Tag 00EJ}]{stacks-project}.

We proceeded in reverse order, with CH proving the ring-theoretic lemma first. Let us state it here:

\begin{lemma}[{\tt Tag 00EJ}]\label{00EJ} Let $R$ be a ring, and say $f_1$, $f_2,\ldots,f_n\in R$ generate the unit ideal. Then the following sequence is exact:
  $$0\longrightarrow R\stackrel{\alpha}\longrightarrow\bigoplus_iR[1/f_i]\stackrel{\beta}\longrightarrow\bigoplus_{i,j}R[1/f_if_j].$$
\end{lemma}

Here the map $\alpha$ is the obvious one, and the map $\beta$ sends $r_i$ to the element whose $(i,j)$ component is $r_i-r_j$.

It was only after we applied this to make progress with Step 1 that we understood the subtleties in what was left. The ``canonical'' identification of basic opens in $\Spec(R[1/f])$ with basic opens in $\Spec(R)$ involved, when identifying global sections, an identification of $R[1/f][1/g]$ with $R[1/fg]$. Of course these rings are canonically isomorphic, but they are not equal. In short, we had a proof of exactness of
$$0\to R[1/f]\to\bigoplus_i(R[1/f])[1/f_i]\to\bigoplus_{i,j}(R[1/f])[1/f_if_j]$$
and we needed\footnote{Strictly speaking we needed a little more than this; in our situation $f_i\in R[1/f]$ rather than $R$, however this special case (a cover of $U$ coming from a cover of $X$) is already enough to illustrate the problem we had.} a proof of exactness of
$$0\to R[1/f]\to\bigoplus_iR[1/ff_i]\to\bigoplus_{i,j}R[1/ff_if_j].$$
To a mathematician, essentially nothing needs to be done here. A mathematician might say ``the diagrams are the same; one is exact, so the other is'' and this would be an acceptable proof. Pressed for more details, a mathematician might offer the following explanation:
$$\xymatrix
{
  0\ar[r]\ar[d]^{=}& R[1/f]\ar[r]\ar[d]^{=}&\bigoplus_i(R[1/f])[1/f_i]\ar[r]\ar[d]^{=}&\bigoplus_{i,j}(R[1/f])[1/f_if_j]\ar[d]^{=}\\
0\ar[r] &R[1/f]\ar[r]&\bigoplus_iR[1/ff_i]\ar[r]&\bigoplus_{i,j}R[1/ff_if_j]
}$$
Here all the vertical maps are canonical and all the horizontal maps are defined in a natural way and hence the squares will all obviously commute.

However in Lean this needs to be checked! Lean has no concept of what it means for an isomorphism to be \emph{canonical} (and looking at \href{https://en.wikipedia.org/wiki/Canonical_map}{the Wikipedia page on canonical maps} one indeed discovers that there seems to be no formal definition of the word), and we needed to explicitly check that the squares in the above diagram commute.

So it came to pass that one line in the Stacks project tag {\tt 01HR} (``Thus we may apply Lemma 10.22.2...We conclude that the sequence is exact'') became several hundred lines of Lean code. In retrospect it would have been much easier to do a diagram chase. However we had a belief that ``everything would follow immediately from the universal property'' and instead went down this route, which was more troublesome than one might expect, because for example the homomorphism $\beta$ above is not a ring homomorphism and hence the universal property cannot be used! We ultimately resorted to proving a lemma showing that there was at most one $R$-algebra map from $R[1/S]$ to $R[1/T]$ and used this to finish. We proved that the squares commuted, and deduced that the presheaf of rings on $\Spec(R)$ was a sheaf. Our goal of proving that affine schemes were schemes was now in sight.

\subsection{Interlude: the univalence axiom.}

Lean does not have the univalence axiom; indeed this axiom contradicts the axiom in Lean saying that all proofs of a theorem of the form $A = B$ are equal. The importance of equality in mathematics is that if $A=B$, and if $C(A)$ is any theorem statement or proof about $A$, then because $A=B$ we have access to the corresponding $C(B)$ about $B$. This process, in a formal proof system, is called \emph{substitution}, or \emph{rewriting}. Univalent systems have a rich concept of equality, where things can be equal in more than one way, whilst preserving this substitution property. The univalence axiom essentially says that if $A$ and $B$ are structures like rings then there's a bijection between isomorphisms $A\to B$ and proofs of $A=B$. 

In particular, in a univalent system, the ring isomorphism $R[1/fg]\cong R[1/f][1/g]$ can be promoted to an equality (for this richer concept of equality) and now it looks on the face of it that a rewrite would be able to make progress. However, univalence does \emph{not} solve the problem which we encountered here. After the substitution in a univalent system, we would have a proof of exactness of one diagram, and we want to prove exactness of another diagram, and the diagrams now have the same objects, however we need to check that they have the same morphisms! Checking this of course boils down to checking that the squares commute, so the lion's share of the work still needs to be done. In our original definition of a scheme, we did the diagram chase ``manually''. However in our second iteration of the definition, we will introduce ideas which enable us to avoid this rewriting problem completely.

\subsection{Affine schemes are schemes.}

With the proof that $\Spec(R)$ is a ringed space, we can now attempt to prove that it is a scheme. This provided one final surprise in our formalisation, and again Reid Barton explained the way around it. Our cover of $\Spec(R)$ by affines is just the identity map $\iota:\Spec(R)\to\Spec(R)$, and all that remains is to show that the presheaves $\cO_X$ and $\iota^*\cO_X$ on $\Spec(R)$ are isomorphic. This boils down to the following: if $U$ is an open subset of $\Spec(R)$ then we need to produce an isomorphism between $\cO_X(U)$ and $\cO_X(\iota(U))$ which commutes with restriction maps.

Our first attempt to do this was the following. Note that $\iota$ is the identity map. Hence $\iota(U)=U$, and thus $\cO_X(\iota(U))=\cO_X(U)$. Let's define the isomorphism to be the identity map. Checking that the diagrams commute should then be straightforward.

But it was not straightforward; one now has to check that $\rho_{\iota(U)\iota(V)}=\rho_{UV}$, and replacing $\iota(U)$ with $U$ in a naive manner caused Lean to give {\tt motive is not type correct} errors. Ultimately we were able to get this working, but what is going on here? Mathematically there seems to be no issue.

We learnt from Mario Carneiro what the problem was. The fact that $\iota(U)=U$ is a trivial theorem, but, perhaps surprisingly, it is not true by definition. The definition of $\iota(U)$ is that it is the set of $x\in X$ such that there exists $u\in U$ with $u=x$. Hence in particular $x\in\iota(U)\iff x\in U$ (this is true by definition), and hence $\iota(U)=U$, because two sets are equal if and only if they have the same elements (this is the axiom of set extensionality). We have proved that $\iota(U)=U$, but along the way we invoked an axiom of mathematics and in particular the equality is not definitional. This means that rewriting the equality $\iota(U)=U$ can cause technical problems with data-carrying types such as $\cO_X(\iota(U))$ which depends on $\iota(U)$ (although in our case they were surmountable). This is a technical issue with dependent type theory and can sometimes indicate that one is working with the wrong definitions.

Fortunately, in our case, Reid Barton pointed out the following extraordinary trick to us: we can define the map $\cO_X(\iota(U))\to\cO_X(U)$ using \emph{restriction} rather than trying to force it to be the identity! Using restriction means that we need to supply a proof that $U\subseteq\iota(U)$, but this is trivial (and who cares that it uses an axiom, we are not trying to rewrite anything). The fact that the diagram commutes now just boils down to the fact that restriction from $\iota(U)$ to $V$ via $U$ equals restriction from $\iota(U)$ to $V$ via $\iota(V)$, which follows from the presheaf axiom of transitivity of restrictions. This is certainly not the way that one would usually think about this, but it works fine.

\subsection{Conclusions}

The main problem with this first approach was the issue with localisations. The ring $R[1/S]$ was defined as an explicit ring, and theorems were proved about it; later on when it came to apply these theorems, it turned out that in our application we only had a ring isomorphic to $R[1/S]$ rather than our explicit definition on the nose. These issues were solved in our second approach.

\section{The second definition.}

Apparently there's a saying in computer science: ``Build one to throw away''. Having done this, we now knew what we should be doing; we had to develop a better theory of localisation. We now describe how we did this.

\subsection{Localisation}

The error we had initially made was to only define ``the'' localisation $R[1/S]$ of a ring~$R$ at a submonoid~$S$. The localisation is defined up to unique isomorphism, but pinning it down as an explicit set (or more precisely an explicit type, as Lean uses type theory rather than set theory) turned out to be a bad idea. What we needed instead is a predicate {\tt is\_localisation\_by S} on ring homomorphisms $R\to T$, saying that $R\to T$ is isomorphic to $R\to R[1/S]$ in the category of $R$-algebras, or in other words that $T$ is isomorphic to $R[1/S]$ in a manner compatible with the $R$-algebra structure. We will refer to this predicate by saying that the $R$-algebra $T$ is \emph{a localisation of $R$ at~$S$}, as opposed to ``the'' localisation of $R$ at $S$.

AL noted that in fact localising rings at submonoids was not the primitive notion: in fact, as Bourbaki taught us, one should be localising monoids at submonoids, and attaching ring structures later on. AL developed an entire formalised theory of localisation of monoids, with both the ``explicit'' constructions $M[1/S]$ and the ``predicate'' approach, showing that the explicit constructions satisfied the predicate, and proving universal properties both for the explicit construction and the predicate construction. These files now form the foundation of the theory of localisation in {\tt mathlib}.

\subsection{The sheaf on an affine scheme, again}

RFM rewrote schemes from scratch, tidying up the code along the way and also moving away from the disastrous design decision of calling Lean files by their Stacks project tags (tags should be mentioned in docstrings, not filenames; filenames serve a useful organisational purpose). The refactoring to using predicates instead of explicit constructions meant that there was no longer any need to make the substitution of $R[1/f][1/g]$ into a lemma explicitly naming $R[1/fg]$; all one has to do is to prove that $R[1/f][1/g]$ satisfies the predicate for being a localisation of $R$ at the submonoid generated by $fg$, which is straightforward. This shortened the definition by hundreds of lines of code.

The hard work, or so we thought, would be in reproving the ring-theoretic lemma \cite[\href{https://stacks.math.columbia.edu/tag/00EJ}{Tag 00EJ}]{stacks-project} (Lemma~\ref{00EJ}) in this more general form. The issue here is that CH's original formalisation was an explicit computation involving the rings $R[1/f_i]$ and $R[1/f_if_j]$; what we now needed was to \emph{reprove} this lemma for \emph{all} rings satisfying the localisation predicate instead. The beefed-up lemma would then apply in situations where the original lemma did not, meaning that the diagram chase described in section~3.3 would no longer be necessary.

\subsection{The localisation predicate}

Definitions are important when formalising. Here are two definitions of being ``a localisation of~$R$ at~$S$'':

\medskip

{\bf Definition 1.} An $R$-algebra $T$ is \emph{a localisation of~$R$ at~$S$} if there exists an $R$-algebra isomorphism $T\cong R[1/S]$.

\medskip

{\bf Definition 2.} An $R$-algebra $T$ is \emph{a localisation of~$R$ at~$S$} if it satisfies the universal property of $R[1/S]$: if $R\to A$ is any ring homomorphism sending $S$ into $A^\times$, there exists a unique map $T\to A$ making the diagram commute.

\medskip

These definitions are mathematically equivalent, so to a mathematician they may as well both be the definition. However Neil Strickland proposed a third definition:

\medskip

{\bf Definition 3 (Strickland).} An $R$-algebra $f:R\to T$ is a \emph{localisation of~$R$ at~$S$} if it satisfies the following three conditions:

\begin{itemize}
\item $f(S)\subseteq T^\times$;
\item Every $t\in T$ can be written as $f(r)/f(s)$ for some $r\in R$ and $s\in S$;
\item The kernel of $f$ is the annihilator of~$S$.
\end{itemize}

\medskip
It is not difficult to prove that Strickland's definition of the predicate is equivalent to the other two. Indeed an explicit computation using the explicit construction of $R[1/S]$ as a quotient of $R\times S$ shows that it satisfies Strickland's predicate, and conversely any $R$-algebra $T$ satisfying Strickland's predicate admits a map from $R[1/S]$ (by the universal property and the first condition), which can be explicitly checked to be surjective (from the second condition) and injective (from the third condition).

So which definition should one use? This is not a mathematical question, it is what is known as an \emph{implementation issue}, or a design decision. Ultimately of course the goal is to prove that all of the definitions are equivalent. But where does one start? One would like to fix one of them and then develop an \emph{interface}, or \emph{API}, for it. This is a collection of basic theorems about our predicate (some deducing it, some using it) which will ultimately lead to a formal proof that it is equivalent to the other two definitions. In particular, one will sometimes have to \emph{verify} the predicate, and one will sometimes have to \emph{use} it. 

Strickland's predicate had advantages over the other two definitions when it comes to using it. In particular, it does not involve quantification over all rings. Using a property involving quantification over all rings often involves having to construct a ring with nice properties and then applying the universal property to that ring. For example, we invite the reader to prove that any $R$-algebra~$T$ satisfying the universal property of being a localisation (definition~2) also has the property that the kernel of $R\to T$ is the annihilator of~$S$. The shortest way we know how to do this from first principles is to first show that the explicitly constructed ring $R[1/S]$ satisfies the universal property, secondly to prove that two rings satisfying the universal property are isomorphic as $R$-algebras, and finally to do an explicit calculation of the kernel of the map $R\to R[1/S]$ to prove the result. Of course there is no problem formalising this proof, but it gives an idea as to how much work it is to build an API for definition 2.

Ultimately one has to fix one definition of the predicate, and then prove that it is equivalent to the other two; there will always be some work involved. But ``API building'' turned out to be easiest with Strickland's definition, which is the one RFM used.

\subsection{Reproving {\tt 00EJ}}

What remains to be done is to reprove Lemma~\ref{00EJ} not for our rings $R[1/f]$ but more generally for rings satisfying Strickland's predicate. To our surprise, this turned out to be very easy. Indeed, the \emph{only} facts about the $R$-algebras $R[1/S]$ used in the Stacks project proof in 2017\footnote{note now that the proof has been changed; the proof formalised in our project is the one from {\tt algebra.tex} in the Stacks project as it stood in late 2017 after commit {\tt 39135dd}.} were precisely the ones isolated by Strickland! The refactoring was hence far easier than expected. In some sense, what happened here was that we made the statement of the lemma more general, and noted that the same proof still worked.

\subsection{Definition and usage.}

RFM also set up the theory of locally ringed spaces, enabling us to define a scheme as a locally ringed space which was locally affine. The proof that an affine scheme is a scheme now also needs a new proof, namely that the stalks of $\cO_X$ were local rings. This needed some of the theory of filtered colimits, but this did not present any particular problems.

KL used the definitions of the project in this form to do the ``gluing sheaves'' Exercise II.1.22 in Hartshorne's algebraic geometry textbook~\cite{hartshorne} and also to prove that $\Spec$ was adjoint to the global sections functor on locally ringed spaces \cite[\href{https://stacks.math.columbia.edu/tag/01I1}{Tag 01I1}]{stacks-project}. However when it came to start thinking about porting our definition into Lean's mathematics library, we ran into a problem. Part of our code was still sub-optimal: we had a definition of a sheaf of types, and a definition of a sheaf of rings. These notions should be unified under some more general notion of a sheaf of objects in a category. Thus the third definition was born.

\section{The third definition: sheaves done correctly.}

\subsection{Category theory in {\tt mathlib}}

SM has been the main driving force behind a gigantic category theory library in Lean. In early 2018, a combination of his work not yet being mature enough, and there not being much infrastructure to support easily working on branches of {\tt mathlib}, meant that we had developed our own definitions of sheaves rather than using category theory. By 2020 this was no longer the case, and SM's definition of a presheaf taking values in a category was the natural thing to use for the ``official'' definition of a scheme. We had seen with our own eyes the problems of having one definition of presheaves of types and another definition of presheaves of rings -- we were constantly having to prove results for presheaves of types and then prove them again for presheaves of rings. Ultimately we had to completely refactor the sheaf part of the story, but we took this opportunity to introduce category theory more generally. Another advantage of the third definition is that rather than working in a new project with {\tt mathlib} as a dependency, we worked directly on a branch of {\tt mathlib}, and ultimately the code ended up as part of {\tt mathlib} meaning that it will not quickly rot and die, as is very very common with Lean code which is not in {\tt mathlib} ({\tt mathlib} is not backwards compatible -- it is still in some sense an experimental library and occasionally big design decisions are changed in an attempt to make it better).

\subsection{Changes made to the definition.} We go through the final definition, pointing out how it differs from the previous version. In this definition, a scheme is an object in the category {\tt LocallyRingedSpace} which is locally isomorphic to an affine scheme. The abstract nonsense presented no surprises; the work left, when defining a scheme and proving that affine schemes are schemes, is to define the sheaf of rings on $\Spec(R)$ in this language. This time it was decided to not follow the Stacks project approach via sheaves on a basis, but to instead define the sheaf directly following~\cite{hartshorne}, where an element of $\cO_X(U)$ is a dependent function taking $u\in U$ to an element of $R_u$, the localisation of $R$ at the prime ideal~$u$, subject to the condition that locally the function can be written as $r/s$ with $r,s\in R$ and $s$ not vanishing near $u$. The advantage of this approach is that it is clear that $\cO_X$ is a sheaf; the ring-theoretic Lemma {\tt 00EJ} is not used at all in this approach. Indeed, in this set-up, {\tt 00EJ} is used to prove $\cO_X(\Spec(R))=R$.

\subsection{Sheaves and categories}

SM defined the notion of a sheaf on a topological space taking values in any category which has products.
The sheaf condition on a presheaf is that the usual ``sheaf condition diagram'' is an equalizer.
% https://github.com/leanprover-community/mathlib/blob/25d83437b19a118aa47db6ead941a4edcee1481a/src/topology/sheaves/sheaf.lean#L68
It is then a theorem that a sheaf of rings is a presheaf of rings whose underlying presheaf of types is a sheaf of types. 
% https://github.com/leanprover-community/mathlib/blob/25d83437b19a118aa47db6ead941a4edcee1481a/src/topology/sheaves/forget.lean#L128
Since then, Bhavik Mehta has defined the concept of a sheaf of types on an arbitrary site, 
% https://github.com/leanprover-community/mathlib/blob/84f99388bec302209eaf0dc70bfd1358d2a0ff46/src/category_theory/sites/sheaf.lean#L604
and more generally the notion of a sheaf of objects of an arbitrary category on a site.
% https://github.com/leanprover-community/mathlib/blob/master/src/category_theory/sites/sheaf.lean
All of this work is now in {\tt mathlib}. These definitions open the door to defining and proving theorems about the \'etale cohomology of schemes in Lean.

\subsection{A word on {\tt mathlib}}

We mentioned above that this definition of a scheme made it into Lean's mathematics library {\tt mathlib}, which means that all the code in it was subject to scrutiny by the library maintainers. A reasonable analogy would be that {\tt mathlib} is like a journal, and the maintainers are like the editors. After three iterations of the definition of a scheme, it was in good enough shape for a ``pull request'' to be made to {\tt mathlib}. (A pull request is the open source project equivalent of submitting an article to an editor!) Note that the multiple iteration procedure discussed in this paper is not usual -- but we were beginners in 2017 with very little information to guide us on how to make schemes on a computer, and KB's poor initial design decisions were because of this. Nowadays we have a much better understanding of how to put modern mathematics into Lean's dependent type theory.

The advantage of getting the definition into {\tt mathlib} is that now it is guaranteed to compile for the duration of Lean 3's lifespan, because if someone makes changes to the library which break it, it will be down to that same someone to fix it. This is in stark contrast to the first two definitions, which compile with very old versions of Lean and {\tt mathlib}, and would almost certainly not compile with modern versions.

\section{Conclusion and acknowledgements.}

We formalised the definition of a scheme, and along the way learnt a lot about how to formalise localisation in dependent type theory. The definition is now in Lean's mathematics library {\tt mathlib} (see \href{https://github.com/leanprover-community/mathlib/blob/b77916dbb7bf5ad0c769b642a4dfde0f2531e49c/src/algebraic_geometry/Scheme.lean#L36-L38}{here}).

We are very grateful to Mario Carneiro, Reid Barton and Neil Strickland, all of whom provided major insights to help us on our way.

\bibliographystyle{amsalpha}
\bibliography{fsbib}
\end{document}